\documentstyle[11pt,psfig,amssymb,amsmath]{article}

\font\smallit=cmti10

\newcommand{\be}{\begin{equation}}
\newcommand{\ee}{\end{equation}}
\newcommand{\bea}{\begin{eqnarray}}
\newcommand{\eea}{\end{eqnarray}}

\def\1#1{^{(#1)}}

\begin{document}
\begin{center}
{\Large\bf On the P\'olya Enumeration Theorem}
\vskip 10pt
{\bf Leonid G. Fel and Yoram Zimmels}\\
{\smallit Department of Civil and Environmental Engineering,
Technion, Haifa 32000, Israel}\\
{\tt lfel@techunix.technion.ac.il}\\
\end{center}

\def\be{\begin{equation}}
\def\ee{\end{equation}}  
\def\bea{\begin{eqnarray}}   
\def\eea{\end{eqnarray}}
\def\p{\prime}

\begin{abstract}
Simple formulas for the number of different cyclic and dihedral necklaces 
containing $n_j$ beads of the $j$-th color, $j\leq m$ and $\sum_{j=1}^mn_j=N$, 
are derived.
\end{abstract}


Among a vast number of counting problems one of the most popular is {\em a 
necklace} enumeration. A {\em cyclic necklace} is a coloring in $m$ colors of 
the vertices of a regular $N$--gon, where two colorings are equivalent if one 
can be obtained from the other by a cyclic symmetry $C_N$, e.g. colored beads 
are placed on a circle, and the circle may be rotated (without reflections). A 
basic enumeration problem is then: for given $m$ and $N=\sum_{j=1}^mn_j$, how 
many different cyclic necklaces containing $n_j$ beads of the $j$-th color 
are there. The answer follows by an application of the P\'olya's theorem 
\cite{har94}: the number $\gamma(C_N,{\bf n}^m)$ of different cyclic 
necklaces is the coefficient of $x_1^{n_1}\cdot\ldots\cdot x_m^{n_m}$ in the 
cycle index 
\begin{eqnarray}
Z_{C_N}(x_i)=\frac{1}{N}\sum_{g|N}\phi(g)X_g^{N/g}\;,\;\;\;\;
X_g=x_1^g+\ldots+x_m^g\;,\label{rip1}
\end{eqnarray}
where $\phi(g)$ denotes the Euler totient function and ${\bf n}^m$ denotes a 
tuple $(n_1,\ldots,n_m)$.

In this article we prove that
\begin{eqnarray}
\gamma\left(C_N,{\bf n}^m\right)=\frac{1}{N}\sum_{d|\Delta}\phi(d)P\left({\bf 
k}^m\right)\;,\;\;\;\mbox{where}\;\;\;P\left({\bf k}^m\right)=\frac{\left(k_1+
\ldots+k_m\right)!}{\prod_{j=1}^mk_j!}\;,\;\;\;k_j=\frac{n_j}{d}\;,\label{rip2}
\end{eqnarray}
and $\Delta$ denotes a great common divisor $\gcd {\bf n}^m$ of the tuple 
${\bf n}^m$. We denote also ${\bf k}^m=(k_1,\ldots,k_m)$.

Note that the term $x_1^{n_1}\cdot\ldots\cdot x_m^{n_m}$ does appear only once 
in {\em the multinomial series expansion} (MSE) of (\ref{rip1}) with a weight 
$P\left({\bf n}^m\right)$ when $g=1$, 
\begin{eqnarray}
X_1^N\;\;\longrightarrow\;\;P\left({\bf n}^m\right)x_1^{n_1}\cdot\ldots\cdot 
x_m^{n_m}\;,\;\;\;\mbox{where}\;\;\;N=n_1+\ldots+n_m\;.\label{rip3}
\end{eqnarray}
Show that for $g>1$ the polynomial $Z_{C_N}\left(x_i\right)$ contributes in 
$\gamma\left(C_N,{\bf n}^m\right)$ if and only if $\Delta>1$. We prove that if 
$g|N$ and $g\not |\Delta$ then the term $x_1^{n_1}\cdot\ldots\cdot x_m^{n_m}$ 
does not appear in MSE of (\ref{rip1}). 

\noindent
Denote $N/g=L$, $1<L<N$ and consider MSE of (\ref{rip1})
\begin{eqnarray}
X_g^L=\sum_{l_i\geq 0}^{l_1+\dots + l_m=L}
P\left({\bf l}^m\right)x_1^{gl_1}\cdot\ldots\cdot x_m^{gl_m}\;,\label{rip4}
\end{eqnarray}
where ${\bf l}^m$ denotes a tuple $(l_1,\ldots,l_m)$. However MSE in 
(\ref{rip4}) does not contribute in $\gamma\left(C_N,{\bf n}^m\right)$ since 
$g\not |\Delta$, i.e. we cannot provide such $g$ that $gl_i=n_i$ holds for 
all $i=1,\dots,m$. Thus, we have reduced expression (\ref{rip1}) by summing 
only over the divisors $d$ of $\Delta$,
\begin{eqnarray}
Z_{C_N}\left(x_i\right)=\frac{1}{N}\sum_{d|\Delta}\phi(d)X_d^{N/d}\;.
\label{rip5}
\end{eqnarray}
Denoting $k_j=n_j/d$, $N/d=K=k_1+\ldots+k_m$, and considering MSE of 
(\ref{rip5}) we obtain
\begin{eqnarray}
X_d^K\;\;\longrightarrow\;\;P\left({\bf k}^m\right)x_1^{dk_1}\cdot\ldots\cdot 
x_m^{dk_m}=P\left({\bf k}^m\right)x_1^{n_1}\cdot\ldots\cdot x_m^{n_m}\;.
\label{rip6}
\end{eqnarray}
Combining (\ref{rip5}) and (\ref{rip6}) we arrive at (\ref{rip2}).

It is easy to extend the explicit formula (\ref{rip2}) to the case of {\em 
dihedral necklaces} where two colorings are equivalent if one can be obtained 
from the other by a dihedral symmetry $D_N$, e.g. colored beads are placed on
a circle, and the circle may be rotated and reflected. Start with the cycle 
indices \cite{har94}
\begin{eqnarray}
2Z_{D_N}\left(x_i\right)=Z_{C_N}(x_i)+\left\{\begin{array}{c}X_1X_2^L\;,\;\;
\mbox{if}\;\;N=2L+1\;,\\\frac{1}{2}\left(X_1^2X_2^{L-1}+X_2^L\right)\;,\;\;
\mbox{if}\;\;N=2L\;.\end{array}\right.\label{rip7}
\end{eqnarray}
If $N=2L+1$ we have to distinguish two different cases.
\begin{enumerate}
\item There is one odd integer $n_j=2a_j+1\in {\bf n}^m$, while the rest of 
$n_i$ are even, $n_i=2a_i$, $L=\sum_{i=1}^ma_i$,
\begin{eqnarray}
\gamma\left(D_N,{\bf n}^m\right)=\frac{1}{2}\left[P\left({\bf n}^m\right)+
P\left({\bf a}^m\right)\right]\;,\;\;\;\mbox{where}\;\;\;{\bf a}^m=(a_1,
\ldots,a_j,\ldots,a_m)\;.\label{rip8}
\end{eqnarray}
\item There is more than one odd integer $n_j=2a_j+1\in {\bf n}^m$, $1\leq j
\leq m$,
\begin{eqnarray}
\gamma\left(D_N,{\bf n}^m\right)=\frac{1}{2}\gamma\left(C_N,{\bf n}^m\right)
\;.\label{rip9}
\end{eqnarray}
\end{enumerate}
If $N=2L$ we have to distinguish three different cases.
\begin{enumerate}
\item All integers $n_j\in {\bf n}^m,j<m$ are even, $n_j=2b_j$, and 
$L=\sum_{i=1}^mb_i$, ${\bf b}^m=(b_1,\ldots,b_m)$,
\begin{eqnarray}
\gamma\left(D_N,{\bf n}^m\right)=\frac{1}{2}\gamma\left(C_N,{\bf n}^m\right)+
\frac{1}{4}\sum_{q=1}^mP\left({\bf b}_q^{m}\right)+\frac{1}{4}P\left({\bf 
b}^m\right)\;,\;\;\;\mbox{where}\label{rip10}
\end{eqnarray}
$${\bf b}_1^m=(b_1-1,b_2,b_3,\ldots,b_m)\;,\;\;{\bf b}_2^m=(b_1,b_2-1,b_3,
\ldots,b_m)\;,\;\ldots\;,\;{\bf b}_m^m=(b_1,b_2,b_3,\ldots,b_m-1)\;.
$$
\item There is one pair of odd integers, $n_{j_1,j_2}=2c_{j_1,j_2}+1\in 
{\bf n}^m$, while the rest of $n_i$ are even, $n_i=2c_i$,
\begin{eqnarray}
\gamma\left(D_N,{\bf n}^m\right)=\frac{1}{2}\left[P\left({\bf 
n}^m\right)+P\left({\bf c}^{m}\right)\right]\;,\;\;\;\mbox{where}\;\;\;{\bf 
c}^m=(c_1,\ldots,c_{j_1},\ldots,c_{j_2},\ldots,c_m),\label{rip11}
\end{eqnarray}
and $L=1+c_1+\dots+c_{j_1}+\dots+c_{j_2}+\dots+c_m$.
\item There is more than one pair of odd integers $n_{j_1,j_2}=2c_{j_1,j_2}+1
\in {\bf n}^m$, $1\leq j_1,j_2\leq m$,
\begin{eqnarray}
\gamma\left(D_N,{\bf n}^m\right)=\frac{1}{2}\gamma\left(C_N,{\bf n}^m\right)
\;.\label{rip12}
\end{eqnarray}
\end{enumerate}

\end{document}